\documentclass[a4paper,11pt]{article}

\usepackage{bm}
\usepackage{anysize}
\usepackage{amsthm}
\usepackage{amssymb}
\usepackage{amsmath}
\usepackage{fancybox}
\usepackage{bbm}
\usepackage{amsfonts}
\usepackage[T1]{fontenc}
\usepackage{latexsym}
\usepackage{makeidx}
\usepackage{cite}
\usepackage[numbers]{natbib}
\usepackage{esint}

\usepackage[pdftex]{graphicx}
\usepackage[pdftex]{color}

\usepackage[all]{xy}

\usepackage{stackengine}
\stackMath

\newtheorem{teo}{Theorem}
\newtheorem{lem}{Lemma}

\newcommand*{\m}[1]{\underline{#1}}

\newcommand{\fd}{\rightarrow}

\newcommand{\inc}{\subset}

\newcommand{\ba}{\overline}
\newcommand{\onda}{\widetilde}

\newcommand{\al}{\alpha}

\newcommand{\lan}{\lambda}

\newcommand{\del}{\delta}

\newcommand{\gam}{\gamma}

\newcommand{\Om}{\Omega}

\newcommand{\Z}{\mathbb{Z}}
\newcommand{\N}{\mathbb{N}}

\newcommand{\R}{\mathbb{R}}
\newcommand{\C}{\mathbb{C}}

\newcommand{\I}{\mathbbm{1}}

\newcommand{\pa}{\partial}

\newcommand{\el}{\ell}

\newtheorem{remark}{Remark}[section]

\def\pf{\par\noindent {\em Proof.}~\par\noindent}

\def\lim{\mathop{\mbox{\normalfont lim}}\limits}

\def\tr{\mathop{\mbox{\normalfont tr}}\nolimits}

\def\pf{\par\noindent {\em Proof. }}

\def\pa{\partial}

\begin{document}

\date{}

\title{Representation formulae for the determinant in a neighborhood of the identity}
\small{
\author
{Denis Constales, Al\'i Guzm\'an Ad\'an}
\vskip 1truecm
\date{\small  Clifford Research Group, Department of Electronics and Information Systems, \\ Faculty of Engineering and Architecture, Ghent University, Krijgslaan 281, 9000 Gent, Belgium. \\
{\tt Denis.Constales@UGent.be,}  \hspace{.3cm}  {\tt Ali.GuzmanAdan@UGent.be}}

\maketitle

\begin{abstract} 
We prove an integral representation and a power series expansion for the function $\det(A)^{-1}$ in a small neighborhood of the identity matrix. Both results are closely linked to the formula for the change of coordinates of the Dirac delta distribution in $\mathbb{R}^m$.




\noindent

\vspace{0.3cm}

\small{ }
\noindent
\textbf{Keywords.} Determinant,  Taylor series, complex analysis, Dirac distribution\\
\textbf{Mathematics Subject Classification (2010).} 	15A15, 41A58, 	30E20

\noindent
\textbf{}
\end{abstract}


\section{Introduction}
In this manuscript we prove two representation formulas for the function $\det(A)^{-1}$ in a small neighborhood of the identity matrix. Let us start by  describing our results.

Let $\C^{k\times k}$ the algebra of complex matrices $M=\{m_{r,\ell}\}_{r,\ell=1\ldots, k}$ of order $(k\times k)$ with identity $\I_k$. 
The {Frobenius} norm of $M\in\C^{k\times k}$ is defined as $\|M\|=(\sum_{r,\ell=1}^k |m_{r,\ell}|^2)^{1/2}$, where $|\cdot|$ denotes the Euclidean norm in $\C$. Associated to each row vector $M_r=(m_{r,1},\ldots, m_{r,k})$, we consider a multi-index $\al_r=(\al_{r,1},\ldots, \al_{r,k})\in \N_0^k$ where $\N_0$ denotes the set of non-negative integers. As usual, we denote $|\al_r|=\al_{r,1}+\ldots+ \al_{r,k}$, $\al_r!=\al_{r,1}!\cdots \al_{r,k}!$ and $M_r^{\al_r}=m_{r,1}^{\al_{r,1}}\cdots m_{r,k}^{\al_{r,k}}$. We shall also consider the multi-index sum $\al_1+\ldots+\al_k=(\sum_r \al_{r,1}, \ldots, \sum_r \al_{r,k})$. In general, for milti-indices $I=(i_1, \ldots, i_k)$ and $J=(j_1, \ldots, j_k)$ we have $I\pm J= (i_1\pm j_1, \ldots, i_k\pm j_k)$. We say that $I\leq J$ if $i_r \leq j_r$ for all $r=1,\ldots, k$.

The multi-indices $\al_r$ give rise to the multi-index matrix ${\bf \al}=\{\al_{r,\ell}\}_{r,\ell=1\ldots, k}\in  \N_0^{k \times k}$. We thus denote $|\al|=|\al_{1}|+\ldots+ |\al_{k}|$, $\al!=\al_{1}!\cdots \al_{k}!$ and $M^{\bf \al}=M_{1}^{\al_{1}}\cdots M_{k}^{\al_{k}}$. We also introduce the differential operator 
  \[\pa_M^\al=\prod_{r,\ell=1}^k \pa_{m_{r,\ell}}^{\al_{r,\ell}},\]
 which is the so-called Fischer dual of the monomial $M^\al$.
 
 With the above notation, the Taylor series of the function $\det(\I_k+M)^{-1}$ around the point $M=0$ can be written as
 \begin{equation}\label{Taylor1st}
 \frac{1}{\det(\I_k+M)}= \sum_{{\bf \al}\in \N_0^{k \times k}} \frac{M^{\bf \al}}{{\bf \al}!} \; \pa_M^\al \left[\frac{1}{\det(\I_k+M)}\right]\Bigg|_{M=0}.
 \end{equation}
Our first goal is to {explicitly write down the above formula}, i.e.\ to compute the derivatives $\pa_M^\al \left[\frac{1}{\det(\I_k+M)}\right]\Big|_{M=0}$. In particular, we will prove the following result.
\begin{teo}\label{Teo1}
Let $M\in \C^{k\times k}$ be such that $\|M\|\leq 1/k$. Then the Taylor series (\ref{Taylor1st}) {converges and} has the form 
 \begin{equation}\label{For1}
 \frac{1}{\det(\I_k+M)}=  \sum_{J\in \N_0^{k }} (-1)^{|J|} J! \sum_{{\substack{\al_1+\ldots+\al_k=J \\ |\al_r|=j_r}}} \frac{M_{1}^{\al_{1}}\cdots M_{k}^{\al_{k}}}{\al_{1}!\cdots \al_{k}!},
 \end{equation}
 where the above sum runs over all multi-indices $J=(j_1, \ldots, j_k)\in  \N_0^{k }$ and all multi-index matrices ${\bf \al}$ such that $\al_{1,r}+\cdots+ \al_{k,r}=|\al_r|=j_r$ for all $r=1,\ldots, k$. This is, 
 ${\bf \al}\in \N_0^{k \times k}$  {is} such that the sum of its $r$-th row {equals the sum} of its $r$-th column for all $r=1,\ldots, k$. 
 \end{teo}
 This constitutes a generalization to higher dimensions of the convergence of the geometric series 
 \[\sum_{j=0}^\infty (-1)^j z^j=\frac{1}{1-z}, \;\;\; \mbox{ when } \;\;\; z\in \C, \; |z|<1.\]
 {It also provides a detailed expression for the full expansion of the well-known formula (see e.g. \cite{MR1231339}) 
 \[ \frac{1}{\det(\I_k+M)} = \exp\left(-\ln \left(\det (\I_k+M)\right)\right) = \exp\left( \sum_{j=0}^\infty (-1)^j\,\frac{\tr(M^j)}{j}\right),\]
 where $\tr(A)$ is the usual matrix trace of $A$. For a detailed account on this and other matrix analysis results we refer the reader (without claiming completeness) to the works \cite{MR1698873, MR1477662, MR3331229}.}
 
 {\noindent Theorem \ref{Teo1} also yields that $\pa_M^\al \left[\frac{1}{\det(\I_k+M)}\right]\Big|_{M=0}$ is different from zero if and only if the sum per row equals the sum per column  in the multi-index matrix $\al$, namely:}
 \begin{equation*}
\pa_M^\al \left[\frac{1}{\det(\I_k+M)}\right]\Bigg|_{M=0} = \begin{cases} (-1)^{|\al|} \, |\al_1|! \cdots |\al_k|! & \mbox{ if } \;\;\al_1+\ldots+\al_k =(|\al_1|, \ldots, |\al_k|), \\ \\
0 & \mbox{otherwise}.\end{cases}
 \end{equation*}
  {In addition, it allows to compute the inverse of the characteristic polynomial $\det(M-\lan \I_k)$ as a power series of $\lan$. Indeed, if $|\lan|\leq (k\|M\|)^{-1}$, one obtains from (\ref{For1}) that 
 \[\frac{1}{\det(M-\lan \I_k)}= 
  \sum_{j=0}^\infty c_j \, \lan^{k-j} , \;\;\; \mbox{ with } \;\;\;  c_j=(-1)^{j+k} \sum_{|J| =j }  J! \sum_{{\substack{\al_1+\ldots+\al_k=J \\ |\al_r|=j_r}}} \frac{M_{1}^{\al_{1}}\cdots M_{k}^{\al_{k}}}{\al_{1}!\cdots \al_{k}!}. \]}
 
{The main motivation for formula (\ref{For1}) comes for the Taylor series expansion of the Dirac distribution, see Section \ref{Sec2}. However, in order to rigorously prove this result, we will need the following integral representation for $\det(\I_k+M)^{-1}$.}
 \begin{teo}\label{Teo2}
Let $A\in  \C^{k\times k}$ be a matrix such that $\|A-\I_k\|<1/k$ and consider the linear transformation $A\m{z}=\m{w}$,  where $\m{z}=(z_1,\ldots, z_k)^T$ and $\m{w}=(w_1,\ldots, w_k)^T$ are vector variables in $\C^k$. Then
 \begin{equation}\label{For2}
\frac{1}{\det(A)}= \frac{1}{(2\pi i)^k} \oint_{(\pa D)^k} \frac{1}{w_1 \cdots w_k}\, dz_1\cdots dz_k,
 \end{equation}
where $\ \oint_{(\pa D)^k}= \oint_{\pa D} \cdots \oint_{\pa D}$ and $\oint_{\pa D}$ denotes the contour integral along the boundary  $\pa D=\{z\in \C: |z|=1\}$ of the unit disk.
\end{teo}
\begin{remark}
As the statements in the  theorems \ref{Teo1} and \ref{Teo2} announce, the condition $\|A-\I_k\|<1/k$ implies that $\det(A)\neq 0$. Indeed, the power series 
$A^{-1}=\left(\I_k+(A-\I_k)\right)^{-1}= \sum_{\ell=0}^\infty (-1)^\el (A-\I_k)^\el$
converges absolutely if $k>1$. Therefore, the inverse matrix $A^{-1}$ exists. For $k=1$, it is clear that $\|A-1\|<1$ implies $A\neq 0$.

{The proximity of $A=\I_k+M$ to $\I_k$ is important for these results to hold. For example, in the class of nilpotent matrices (that are somehow close to the origin) formula (\ref{For1}) does not hold in general. Indeed, if $M$ is nilpotent then $\I_k+M$ is invertible, however the series (\ref{For1}) may not converge. See, for example, {the case where} $M=\left(\hspace{-.1cm}\begin{array}{cc} 0& 1 \\ 0 & 0 \end{array}\hspace{-.1cm}\right)$.}

{The statement of Theorem \ref{Teo2} also announces that none of the $w_j$'s becomes zero as $\m{z}$ varies in the boundary $(\pa D)^k$ of the unit polydisc. This is another reason why the proximity of $A$ to $\I_k$ is crucial. As a matter of fact, $(\pa D)^k$ can be continuously deformed into $A(\pa D)^k$ without ever touching one of the axis $z_j=0$, see our Lemma \ref{Homology}. }
\end{remark}

{This paper is organized as follows. In section \ref{Sec2}, we informally discuss the main motivations for Theorem \ref{Teo1}. This result is closely connected to the formula for the change of coordinates of the Dirac delta distribution.  Most of the heuristic and motivational reasoning  in that section is not completely rigorous. In the later sections, all of our results will be rigorously proved. Our strategy to prove Theorem \ref{Teo1} is to show first that it is equivalent to Theorem \ref{Teo2} and then proceed to prove the latter. In Section \ref{Sec3}, we prove  that  the left-hand side expressions in (\ref{For1}) and (\ref{For2}) coincide when $A=\I_k+M$, completing in this way the first step of our strategy. The second step is completed in Section \ref{Sec4}, where we prove Theorem \ref{Teo2} by means of the method of Gaussian elimination. Finally, in Section \ref{Sec5}, we provide an alternative proof by showing that the manifolds $(\pa D)^k$ and $A(\pa D)^k$ are in the same homology class in the space $(\C\setminus\{0\})^k$ and therefore, integrals of closed differential $k$-forms over these manifolds remain  invariant. At the end, we briefly discuss the connection between Theorem \ref{Teo2} and the formula for the change of coordinates of the delta distribution.}



\section{Motivations and informal discussion}\label{Sec2}
In this section we sketch the main motivation behind formula (\ref{For1}). This does not lead to a rigorous proof but it  shows how this result is linked with the change of coordinates of the Dirac delta distribution. In the next sections, we will rigorously prove theorems \ref{Teo1} and \ref{Teo2}.

Given the real vector variable $\m{x}=(x_1,\ldots, x_k)^T$, we consider the Dirac distribution $\del(\m{x})=\del(x_1)\ldots \del(x_k)$ in the $k$-dimensional Euclidean space $\R^k$. For any non-singular real matrix $\I_k+M\in \textup{GL}(k,\R)$ it is known that (see e.g. \cite{MR0166596})
\begin{equation}\label{delPri}
\del(\m{x}+M \m{x}) = \del\left(\left(\I_k+M\right)\m{x}\right) = \frac{\del(\m{x})}{|\det(\I_k+M)|}.
\end{equation}
On the other hand, if we formally write down the Taylor series expansion of the left-hand side, we get 
\[\del(\m{x}+M \m{x})  =  \sum_{J\in \N_0^{k }} \frac{\left(M \m{x}\right)^{J}}{{J}!} \; \del^{(J)}(\m{x}),\]
where $\del^{(J)}(\m{x})= \del^{(j_1)}(x_1)\ldots \del^{(j_k)}(x_k)$ for the multi-index $J=(j_1, \ldots, j_k)$. We recall that $\left(M \m{x}\right)^{J} = \left(M \m{x}\right)_1^{j_1} \ldots \left(M \m{x}\right)_k^{j_k}$ where $\left(M \m{x}\right)_r=\sum_{\el=1}^k m_{r,\el} x_\el$ is the $r$-th component of the vector $M\m{x}$. By the multinomial theorem we have 
\[\frac{\left(M \m{x}\right)_r^{j_r}}{j_r!} = \sum_{|\al_r|=j_r} \frac{M_r^{\al_r} \, \m{x}^{\al_r}}{\al_r!}, \;\;\; \mbox{ with } \;\;\; {\m{x}}^{\al_r}=x_{1}^{\al_{r,1}}\cdots x_{k}^{\al_{r,k}},\]
and therefore
\begin{equation}\label{DelMain}
\del(\m{x}+M \m{x})  =  \sum_{J\in \N_0^{k}} \left( \sum_{|\al_r|=j_r} \frac{M_{1}^{\al_{1}}\cdots M_{k}^{\al_{k}}}{\al_{1}!\cdots \al_{k}!} {\m{x}}^{\al_1+\ldots \al_k} \right) \; \del^{(J)}(\m{x}).
\end{equation}
Let us consider $I=\al_1+\cdots \al_k\in \N_0^k$, it is a known result that (see e.g. \cite{MR0166596})
\[
\m{x}^I \, \del^{J}=\begin{cases} (-1)^I \frac{J!}{(J-I)!}  \del^{(J-I)}(\m{x}), & \mbox{ if } I\leq J, \\ 0, &  \mbox{ otherwise.} \end{cases}
\]
But our multi-index $I$ satisfies $|I|=|J|$. Thus, in this case, the condition $I\leq J$ implies $I=J$. Therefore, formula (\ref{DelMain}) can be rewritten as 
\[\del(\m{x}+M \m{x})  = \left(\sum_{J\in \N_0^{k }} (-1)^{|J|} J! \sum_{{\substack{\al_1+\ldots+\al_k=J \\ |\al_r|=j_r}}} \frac{M_{1}^{\al_{1}}\cdots M_{k}^{\al_{k}}}{\al_{1}!\cdots \al_{k}!}\right) \del(\m{x}). \]
Comparing this with (\ref{delPri}), it follows that formula (\ref{For1}) should hold for some suitable class of {matrices} $M$. The statement of Theorem \ref{Teo1} is stronger than this guess. In fact, it explicitly describes a class of matrices for which this formula holds and it states the result for complex matrices. In Section \ref{Sec5}, we shall make the link between the other representation formula (\ref{For2}) and the change of coordinates in the Dirac distribution.


\section{An intermediate step}\label{Sec3}
Before rigurously proving Theorem \ref{Teo1}, we will show that the statements in the theorems \ref{Teo1} and \ref{Teo2} are equivalent. Let $M\in\C^{k\times k}$ be as in Theorem \ref{Teo1} and let us denote  the sum in the right-hand side of (\ref{For1}) by
\[R(M)= \sum_{J\in \N_0^{k }} (-1)^{|J|} J! \sum_{{\substack{\al_1+\ldots+\al_k=J \\ |\al_r|=j_r}}} \frac{M_{1}^{\al_{1}}\cdots M_{k}^{\al_{k}}}{\al_{1}!\cdots \al_{k}!}.\]
Now consider the following {\it "more relaxed"} version of $R(M)$
\begin{equation}\label{For4}
S(M) =  \sum_{J\in \N_0^{k }} (-1)^{|J|} J! \sum_{|\al_r|=j_r} \frac{M_{1}^{\al_{1}}\cdots M_{k}^{\al_{k}}}{\al_{1}!\cdots \al_{k}!}. 
\end{equation}
While the sum $S(M)$ runs over all multi-index matrices ${\bf \al}\in \N_0^{k \times k}$, the sum $R(M)$ considers only those  ${\bf \al}\in \N_0^{k \times k}$ for which the sum of its $r$-th row equals the sum of its $r$-th column for all $r=1,\ldots, k$. Using the multinomial theorem, we can write $S(M)$ in terms of the following geometric series,
\begin{equation}\label{For5}
S(M)=  \sum_{J\in \N_0^{k }}\,  \prod_{r=1}^k \left(-\sum_{\ell=1}^k m_{r,\ell}\right)^{j_r} = \prod_{r=1}^k \, \left(\sum_{j_r=0}^\infty \left(-\sum_{\ell=1}^k m_{r,\ell}\right)^{j_r}\right)= \prod_{r=1}^k \frac{1}{1+ \sum_{\ell=1}^k m_{r,\ell}}. 
\end{equation}
We recall that the power series $\sum_{j_r=0}^\infty \left(-\sum_{\ell=1}^k m_{r,\ell}\right)^{j_r}$ converges uniformly to $\left({1+ \sum_{\ell=1}^k m_{r,\ell}}\right)^{-1}$ since $\|M\| \leq \frac{1}{k}$. Indeed, $\left| \sum_{\ell=1}^k m_{r,\ell} \right| \leq \sum_{\ell=1}^k |m_{r,\ell}| < k \frac{1}{k}=1$. This reasoning also shows that the series $R(M)$ converges absolutely for $\|M\| \leq \frac{1}{k}$.

Our strategy is to apply some transformations to the sum $S(M)$ in order to recover $R(M)$. To that end we first consider the transformation $M \mapsto D^{-1}MD$ where $D=\textup{diag}(z_1, \ldots, z_k)$ is a diagonal matrix whose diagonal entries are in the unit circle, i.e.\ $z_1, \ldots, z_k \in \pa D$. The matrix $D^{-1}MD$ is the result of multiplying the $r$-th
 row of $M$ by $z_r^{-1}$ and the $r$-th column by  $z_r$, $r=1,\ldots, k$. Then for every entry this transformation can be written as
\begin{equation}\label{For6}
m_{r,\ell} \mapsto z_r^{-1}z_\ell \, m_{r,\ell} , \hspace{1cm} r,\ell=1,\ldots, k.
\end{equation}
Let us examine how the sum written in (\ref{For4}) reads after this transformation. Observe that 
\[M_r^{\al_r}=m_{r,1}^{\al_{r,1}}\cdots m_{r,k}^{\al_{r,k}} \;\; \mapsto\;\;  z_r^{-j_r} z_1^{\al_{r,1}}\cdots  z_k^{\al_{r,k}} \,M_r^{\al_r},\]
which yields
\[M_{1}^{\al_{1}}\cdots M_{k}^{\al_{k}} \;\; \mapsto\;\;  \left(\prod_{r=1}^k \displaystyle z_r^{\sum_\ell \al_{\el,r}-j_r}\right)M_{1}^{\al_{1}}\cdots M_{k}^{\al_{k}}.\]
 The only summands in  (\ref{For4}) that remain independent of the $z_r$'s are those satisfying $\sum_{\ell=1}^k \al_{\el,r}=j_r$ for all $r=1,\ldots,k$, or equivalently,  $\al_1+\ldots+\al_k=J$. These are exactly the terms that appear in $R(M)$. Hence we can write
 \[S(D^{-1}MD) = R(M) +T(D^{-1}MD),\]
 where $T(D^{-1}MD)$ is a sum of elements of the form $z_1^{\lan1} \cdots z_k^{\lan_k} c$ such that  $c$ is independent of the $z_r$'s and at least one of the powers $\lan_r\in \Z$ is different from zero. 
 
 Using Cauchy's integral theorem we easily find that 
  \[
 \frac{1}{(2\pi i)^k} \oint_{(\pa D)^k} \frac{T(D^{-1}MD)}{z_1 \cdots z_k}\, dz_1\cdots dz_k = 0
 \]
and
 \[
 \frac{1}{(2\pi i)^k} \oint_{(\pa D)^k} \frac{R(M)}{z_1 \cdots z_k}\, dz_1\cdots dz_k = R(M) \frac{1}{(2\pi i)^k} \oint_{(\pa D)^k} \frac{dz_1\cdots dz_k}{z_1 \cdots z_k} = R(M).
 \]
 Thus
 \[
  \frac{1}{(2\pi i)^k} \oint_{(\pa D)^k} \frac{S(D^{-1}MD)}{z_1 \cdots z_k}\, dz_1\cdots dz_k = R(M).
 \]
 On the other hand, using formulas (\ref{For5}) and (\ref{For6}), we get 
 \[S(D^{-1}MD)=\prod_{r=1}^k \frac{1}{1+ \sum_{\ell=1}^k z_r^{-1}z_\ell \,m_{r,\ell}} = \prod_{r=1}^k \frac{z_r}{z_r+ \sum_{\ell=1}^k z_\ell \,m_{r,\ell}}= \frac{z_1 \cdots z_k}{w_1 \cdots w_k}.\]
 Here the vectors $\m{z}=(z_1,\ldots, z_k)^T$ and $\m{w}=(w_1,\ldots, w_k)^T$ are as in Theorem \ref{Teo2}, i.e.\ $A\m{z}=\m{w}$ with $A=\I_k+M$. Finally, combining the last two formulas, we obtain
 \[R(M)= \frac{1}{(2\pi i)^k} \oint_{(\pa D)^k} \frac{1}{w_1 \cdots w_k}\, dz_1\cdots dz_k,\]
 which proves that the left-hand side expressions in (\ref{For1}) and (\ref{For2}) coincide when $A=\I_k+M$.

\section{Proofs of the main theorems}\label{Sec4}
We now proceed to prove Theorem \ref{Teo2}. To that end, we first need the following lemma.
\begin{lem}\label{Lem1}
Let $A=\{a_{r,\ell}\}_{r,\ell=1,\ldots, k}$ be a matrix in $\C^{k\times k}$ such that $\|A-\I_k\|<1/k$ and consider the linear transformation $A\m{z}=\m{w}$,  where $\m{z}=(z_1,\ldots, z_k)^T \in (\pa D)^k$ and $\m{w}=(w_1,\ldots, w_k)^T$. Then 
\[\frac{1}{2\pi i} \oint_{\pa D} \frac{dz_1}{w_r}=\begin{cases} \displaystyle\frac{1}{a_{1,1}} & \mbox{ if } r=1,   \\ 0 & \mbox{ if }  r=2,\ldots, k.\end{cases}\]
\end{lem}
\pf
Let us start by considering the case $r=2,\ldots, k$. If $a_{r,1}=0$, then $w_r$ does not depend on $z_1$ and the integral is automatically zero. If $a_{r,1}\neq 0$, we have
\[\frac{1}{w_r}=\frac{1}{a_{r,1}\left(z_1+\sum_{\el=2}^k \frac{a_{r,\el}}{a_{r,1}}z_\el\right)}.\]
This is a function of $z_1$ with only one singularity, namely $z_1=-\sum_{\el=2}^k \frac{a_{r,\el}}{a_{r,1}}z_\el$, which lies outside of the unit disk $D=\{z_1\in \C: |z_1|\leq 1\}$. Indeed, by the triangular inequality we obtain
\[
\left|\sum_{\el=2}^k \frac{a_{r,\el}}{a_{r,1}}z_\el \right| = \frac{1}{|a_{r,1}|} \left|\sum_{\el=2}^k {a_{r,\el}} z_\el \right| \geq \frac{1}{|a_{r,1}|} \left(|a_{r,r}| - \left|\sum_{\el\neq 1,r} {a_{r,\el}} z_\el \right|  \right) \geq \frac{1}{|a_{r,1}|} \left(|a_{r,r}| - \sum_{\el\neq 1,r} \left|{a_{r,\el}} \right|  \right).
\]
Now we recall that $|a_{r,\el}-\del_{r,\ell}|<\frac{1}{k}$. In particular, this implies that $\frac{1}{|a_{r,1}|}>k$, $|a_{r,r}|>\frac{k-1}{k}$ and $-|{a_{r,\el}}|> -\frac{1}{k}$ ($\el\neq r$). We thus obtain
\[\left|\sum_{\el=2}^k \frac{a_{r,\el}}{a_{r,1}}z_\el \right| > k\left(\frac{k-1}{k} - \frac{k-2}{k}\right)=1.\]
 This means that $\frac{1}{w_r}$ is a holomorphic function inside the unit disk and therefore $\frac{1}{2\pi i} \oint_{\pa D} \frac{dz_1}{w_r}=0$.
 
\noindent In the case where $r=1$, it suffices to show that the isolated singularity $z_1=-\sum_{\el=2}^k \frac{a_{1,\el}}{a_{1,1}}z_\el$ of $\frac{1}{w_1}$ is inside of the unit disk. One easily observes that
\[\left| \sum_{\el=2}^k \frac{a_{1,\el}}{a_{1,1}}z_\el \right|\leq \sum_{\el=2}^k \frac{|a_{1,\el}|}{|a_{1,1}|} < \sum_{\el=2}^k \frac{\frac{1}{k}}{1-\frac{1}{k}}=1.\]
Then, by the residue theorem we obtain $\frac{1}{2\pi i} \oint_{\pa D} \frac{dz_1}{w_1}=\frac{1}{a_{1,1}}$. $\hfill\square$

\vspace{.3cm}
\noindent {\it Proof of Theorem \ref{Teo2}}

\noindent We proceed by induction on $k\in \N$. For $k=1$ we have $w_1=a_{1,1} z_1$ with $|a_{1,1}-1|<1$. It is then clear that $\frac{1}{2\pi i} \oint_{\pa D} \frac{dz_1}{a_{1,1} z_1}=\frac{1}{a_{1,1}}=\frac{1}{\det(A)}$.

\noindent Let us assume that formula (\ref{For2}) is true for $k-1\in \N$, and let us prove that it also holds for  $k$. To that end we first decompose the function $\frac{1}{w_1\cdots w_k}$ as a sum of partial fractions with respect to $z_1$, i.e.\
\[\frac{1}{w_1\cdots w_k}=\frac{\lan_1}{w_1}+\cdots+\frac{\lan_k}{w_k}\]
where $\lan_1,\ldots, \lan_k$ do not depend on $z_1$. From Lemma \ref{Lem1} we obtain
\[\frac{1}{2\pi i} \oint_{\pa D}\frac{dz_1}{w_1\cdots w_k} = \frac{\lan_1}{a_{1,1}}.\]
We recall that $\lan_1$ is the residue of the function $\frac{a_{1,1}}{w_1\cdots w_k}$ at the singularity $z_1=-\sum_{\el=2}^k \frac{a_{1,\el}}{a_{1,1}}z_\el$. Thus $\lan_1$ can be easily computed to be
\[\lan_1=\frac{1}{w_2\cdots w_k}\bigg|_{z_1=-\sum_{\el=2}^k \frac{a_{1,\el}}{a_{1,1}}z_\el}=\frac{1}{\onda{w}_2 \cdots \onda{w}_k}\]
where $\onda{w}_r$ is the value of $w_r$ when substituting $z_1=-\sum_{\el=2}^k \frac{a_{1,\el}}{a_{1,1}}z_\el$, $r=2,\ldots, k$. Further computations yield
\[
\onda{w}_r =  -\frac{a_{r,1}}{a_{1,1}} \sum_{\el=2}^k a_{1,\el} z_\el +\sum_{\el=2}^k a_{r,\el} z_\el = \sum_{\el=2}^k b_{r,\el} z_\el, \;\;\;\; \mbox{where }\;\;\;\;  b_{r,\el}=a_{r,\el} - \frac{a_{r,1}}{a_{1,1}} a_{1,\el}.
\]
In this way we have obtained
\begin{equation}\label{Int1}
 \frac{1}{(2\pi i)^k} \oint_{(\pa D)^k} \frac{1}{w_1 \cdots w_k}\, dz_1\cdots dz_k =  \frac{1}{a_{1,1}} \frac{1}{(2\pi i)^{k-1}} \oint_{(\pa D)^{k-1}} \frac{1}{\onda{w}_2 \cdots \onda{w}_k}\, dz_2\cdots dz_k,
 \end{equation}
where $\left(\hspace{-.1cm}\begin{array}{ccc}  b_{2,2}& \ldots & b_{2,k} \\  \vdots& \ddots &\vdots \\  b_{k,2}& \ldots & b_{k,k}\end{array}\hspace{-.1cm}\right) \hspace{-.1cm} \left(\hspace{-.1cm}\begin{array}{c}  z_2 \\ \vdots \\  z_k\end{array}\hspace{-.1cm}\right) \hspace{-.1cm}= \hspace{-.1cm}\left(\hspace{-.1cm}\begin{array}{c}  \onda{w}_2 \\ \vdots \\  \onda{w}_k\end{array}\hspace{-.1cm}\right)$ and the matrix $B=\{b_{r,\el}\}_{r,\ell=2,\ldots, k}$ satisfies $\|B-\I_{k-1}\|\leq \frac{1}{k-1}$. Indeed, the matrix $B-\I_{k-1}$ has entries $b_{r,\el}-\del_{r,\ell}=(a_{r,\el}-\del_{r,\ell}) - \frac{1}{a_{1,1}} a_{r,1}a_{1,\el}$ with $r,\ell=2,\ldots, k$. If we donote by $C$ the matrix $\{a_{r,1}a_{1,\el}\}_{r,\ell=2,\ldots, k}$, we get
\[\|B-\I_{k-1}\|\leq \|A-\I_{k}\| +  \frac{1}{|a_{1,1}|}  \|C\|.\]
But $\|C\| = \left(\sum_{r,\el=2}^k |a_{r,1}|^2 |a_{1,\el}|^2 \right)^{\frac{1}{2}}= \left(\sum_{r=2}^k |a_{r,1}|^2 \right)^{\frac{1}{2}} \left(\sum_{\el=2}^k |a_{1,\el}|^2 \right)^{\frac{1}{2}} <  \frac{1}{k^2}$ which implies
\[\|B-\I_{k-1}\|\leq \frac{1}{k} +  \frac{k}{k-1} \frac{1}{k^2} = \frac{1}{k-1}.\]
Now, applying our induction hypothesis on (\ref{Int1}),  we obtain
\[ \frac{1}{(2\pi i)^k} \oint_{(\pa D)^k} \frac{1}{w_1 \cdots w_k}\, dz_1\cdots dz_k = \frac{1}{a_{1,1} \det(B)}. \]
Thus it suffices to prove that $a_{1,1} \det(B)=\det(A)$. This easily follows from the Gaussian elimination process. Indeed,  if we add to the $r$-th row in the matrix $A$ the first row multiplied by $- \frac{a_{r,1}}{a_{1,1}}$ ($r=2,\ldots,k$),  we obtain the matrix
\[\left(\hspace{-.1cm}\begin{array}{cccc} a_{1,1} & a_{1,2} & \ldots & a_{1,k} \\ 0 & b_{2,2}& \ldots & b_{2,k}  \\ \vdots&  \vdots& \ddots &\vdots \\ 0&  b_{k,2}& \ldots & b_{k,k}\end{array}\hspace{-.1cm}\right).\]
This matrix has the same determinant as $A$. We thus obtain that $\det(A)=a_{1,1} \det(B)$, which proves the result. $\hfill\square$
 
 

\section{An alternative proof}\label{Sec5}
We now provide an alternative proof for Theorem \ref{Teo2}. In particular, we will prove the following slightly generalized result. 
 \begin{teo}\label{Teo3}
Consider $A\m{z}=\m{w}$ as in Theorem \ref{Teo2} and let $f(\m{z})$ be a $\C$-valued holomorphic function in $\C^k$. Then
 \begin{equation}\label{For8}
\frac{f(0)}{\det(A)}= \frac{1}{(2\pi i)^k} \oint_{(\pa D)^k}  \frac{f(\m{z})}{w_1 \cdots w_k}\, dz_1\cdots dz_k.
 \end{equation}
\end{teo}
We first observe that $A\m{z}=\m{w}$ yields  $dw_1\cdots dw_k = \det(A) \, dz_1\cdots dz_k$. Effectuating this change of variables in the right-hand side of (\ref{For8}) we obtain
\begin{equation}\label{For9}
 \frac{1}{(2\pi i)^k} \oint_{(\pa D)^k}  \frac{f(\m{z})}{w_1 \cdots w_k}\, dz_1\cdots dz_k =  \frac{1}{\det(A)} \frac{1}{(2\pi i)^k} \oint_{A(\pa D)^k}  \frac{f(A^{-1}\m{w})}{w_1 \cdots w_k}\, dw_1\cdots dw_k.
 \end{equation}
Thus it suffices to prove that $\frac{1}{(2\pi i)^k} \oint_{A(\pa D)^k}  \frac{f(A^{-1}\m{w})}{w_1 \cdots w_k}\, dw_1\cdots dw_k=f(0)$. The idea of this proof is to note that $\frac{f(A^{-1}\m{w})}{w_1 \cdots w_k}\, dw_1\cdots dw_k$ is a closed differential form on  $(\C\setminus \{0\})^k$. Therefore, by Stoke's theorem, the above integral remains invariant when taken over any other manifold in the same homology class of $(\pa D)^k$ in $(\C\setminus \{0\})^k$, { see e.g. \cite{MR2723362}}. This is summarized in the following two lemmas.
\begin{lem}\label{Homology}
Let $A\in  \C^{k\times k}$ be a matrix such that $\|A-\I_k\|<1/k$. Then $(\pa D)^k$ and $A (\pa D)^k$ belong to the same homology class in $(\C\setminus \{0\})^k$, i.e.\ there exists a continuous map $A(t):[0,1]\fd \C^{k\times k}$ such that $A(0)=\I_k$, $A(1)=A$ and $A(t) (\pa D)^k \inc (\C\setminus \{0\})^k$ for every $t\in [0,1]$.
\end{lem}
\pf
Let us consider the map $A(t)=\I_k+t(A-\I_k)$ and $\m{w}(t)=A(t) \m{z}$ with $\m{z}\in (\pa D)^k$ arbitrary. We need to prove that $\m{w}(t)\in  (\C\setminus \{0\})^k$ for every $t\in [0,1]$, or equivalently, that every entry $w_r(t)$ of $\m{w}(t)$  ($r=1,\ldots, k$) is different from zero.

\noindent We recall that $w_r(t)=z_r+t\sum_{j=1}^k (a_{r,j}-\del_{r,j})z_j$. Then by the triangular inequality we have
\begin{align*}
1&= |z_r| = \left|w_r(t) -  t\sum_{j=1}^k (a_{r,j}-\del_{r,j})z_j \right| \leq |w_r(t)| + t \sum_{j=1}^k |a_{r,j}-\del_{r,j}| < |w_r(t)| + t \sum_{j=1}^k \frac{1}{k}.
\end{align*}
Thus $|w_r(t)|>1-t\geq 0$ and therefore $w_r(t)\neq 0$.
$\hfill\square$
 
 \begin{lem}
Let $\gam(\m{w})$ be a $\C$-valued holomorphic function in an open region $\Om \inc \C^k$. Then the differential form $\gam  dw_1\cdots dw_k$ is closed in $\Om$, i.e. $d(\gam  dw_1\cdots dw_k)=0$ where $d$ is the exterior derivative.
 \end{lem}
\pf It is easily seen that the exterior derivative can be written as $d=\pa+\ba{\pa}$ where $\pa= \sum_{j=1}^n \pa_{w_j} dw_j$ and $\ba{\pa}= \sum_{j=1}^n \pa_{\ba{w}_j} d\ba{w}_j$ are given in terms of the classical Cauchy-Riemann operators $ \pa_{w_j}$,  $\pa_{\ba{w}_j}$ and the complex differentials $dw_j$, $d\ba{w}_j$. If we write $w_j=x_j+iy_j$ (with $x_j,y_j$ being real variables) then 
\begin{align*}
\pa_{w_j}&= \frac{1}{2}(\pa_{x_j}-i\pa_{y_j}),   & dw_j &=d{x_j}+id{y_j},   &   \pa_{\ba{w}_j}&= \frac{1}{2}(\pa_{x_j}+i\pa_{y_j}),   & d\ba{w}_j &=d{x_j}-id{y_j}.
\end{align*}
It is clear that $\pa (\gam  dw_1\cdots dw_k)=0$ and, since $\gam$ is holomorphic, we also have that $\ba{\pa} (\gam  dw_1\cdots dw_k)=\sum_{j=1}^n \pa_{\ba{w}_j}[\gam] d\ba{w}_j  dw_1\cdots dw_k =0 $. Therefore $d(\gam  dw_1\cdots dw_k)=0$. 
$\hfill\square$

Using the previous lemma we easily observe that  $\frac{f(A^{-1}\m{w})}{w_1 \cdots w_k}\, dw_1\cdots dw_k$ is a closed differential form on $(\C\setminus \{0\})^k$. This means that its integrals over the homologous manifolds $(\pa D)^k$ and $A(\pa D)^k$ are equal. Going back to formula (\ref{For9}), we finally obtain from Cauchy's theorem that
\begin{align*}
 \frac{1}{(2\pi i)^k} \oint_{(\pa D)^k}  \frac{f(\m{z})}{w_1 \cdots w_k}\, dz_1\cdots dz_k &=  \frac{1}{\det(A)} \frac{1}{(2\pi i)^k} \oint_{A(\pa D)^k}  \frac{f(A^{-1}\m{w})}{w_1 \cdots w_k}\, dw_1\cdots dw_k \\
 &= \frac{1}{\det(A)} \frac{1}{(2\pi i)^k} \oint_{(\pa D)^k}  \frac{f(A^{-1}\m{w})}{w_1 \cdots w_k}\, dw_1\cdots dw_k \\
 &= \frac{f(0)}{\det(A)}.
\end{align*}

\subsection{Connection with the Dirac distribution}
Let us consider the $2k$-dimensional real vector variables $\m{x}=(x_1,\ldots, x_{2k})^T$ and $\m{y}=(y_1,\ldots, y_{2k})^T$ associated to the complex vector variables $\m{z}$ and $\m{w}$ in $\C^k$ by means of 
\[\m{z}=(x_1+ix_{k+1},\ldots, x_k+ix_{2k})^T, \;\;\; \mbox{ and } \;\;\; \m{w}=(y_1+iy_{k+1},\ldots, y_k+iy_{2k})^T,\]
respectively. Equivalently we may write $\m{z}=P \m{x}$ and $\m{w}=P \m{y}$ where $P=\left(\I_k| i\I_k\right)\in \C^{k \times 2k}$. Associated to any complex-linear transformation $\m{w}=A \m{z}$, one finds a real-linear transformation $\m{y}=\Psi(A) \m{x}$, where $\Psi:\C^{k \times k} \fd \R^{2k \times 2k}$ is an algebra morphism given by
\[
\Psi(A_1+iA_2)=  \left(\hspace{-.1cm}\begin{array}{cc} A_1& -A_2 \\ A_2 & A_1 \end{array}\hspace{-.1cm}\right), \;\;\;\;\; A_1,A_2\in \R^{k \times k}.
\]
The determinants of the matrices $A$ and $\Psi(A)$ are linked by the relation $\det(\Psi(A))=|\det(A)|^2$. Indeed, if one considers the matrices $D= \left(\hspace{-.1cm}\begin{array}{cc} \I_k& i\I_k  \\ 0 & \I_k \end{array}\hspace{-.1cm}\right)$ and its inverse $D^{-1}= \left(\hspace{-.1cm}\begin{array}{cc} \I_k& -i\I_k  \\ 0 & \I_k \end{array}\hspace{-.1cm}\right)$, one obtains 
\[D \, \Psi(A) \,D^{-1} = \left(\hspace{-.1cm}\begin{array}{cc} A_1+iA_2 & 0 \\ A_2 & A_1-iA_2 \end{array}\hspace{-.1cm}\right),\]
and therefore $\det(\Psi(A))=\det(D \, \Psi(A) \,D^{-1}) = \det(A_1+iA_2) \det(A_1-iA_2) = |\det(A)|^2$.

Let us define $\del(\m{z})=\del(\m{x})=\del(x_1)\ldots \del(x_{2k})$ and $\del(\m{w})=\del(\m{y})=\del(y_1)\ldots \del(y_{2k})$. Then (see e.g. \cite{MR0166596})
\begin{equation}\label{delCompx}
\del(\m{w}) = \frac{\del(\m{z})}{\det(\Psi(A))}=\frac{\del(\m{z})}{ |\det(A)|^2}.
\end{equation}
In Section \ref{Sec2}, we showed how this formula is linked with Theorem \ref{Teo1}. In this section, we shall make the relation between formula (\ref{delCompx}) and Theorem \ref{Teo3} explicit.

From Theorem \ref{Teo3}, we have for every holomorphic function $f(\m{z})$ that
\begin{equation}\label{delInt}
\left\langle \frac{\del(\m{z})}{\det(A)}, f \right\rangle=\frac{f(0)}{\det(A)}=\frac{1}{\det(A)} \frac{1}{(2\pi i)^k} \oint_{(\pa D)^k}  \frac{f(\m{z})}{w_1 \cdots w_k}\, dw_1\cdots dw_k.
\end{equation}
We now recall that Green's theorem can be written, in terms of the Cauchy-Riemann $\pa_{\ba{z}}=\frac{1}{2}(\pa_x+i\pa_y)$ operator of the complex variable $z=x+iy$, as
\[\oint_{\pa D} g \, dz = 2i \iint_D \pa_{\ba{z}}[g] \, dx dy,\]
where $g$ is a differentiable function in a neighborhood of the unit disc $D$. Applying Green's theorem in each variable $w_j$ in (\ref{delInt}), we obtain
\[
\left\langle \frac{\del(\m{z})}{\det(A)}, f \right\rangle=\frac{1}{\pi^k \det(A)} \int_{D^k}  \pa_{\ba{w}_1} \cdots \pa_{\ba{w}_{2k}}\left[\frac{1}{w_1 \cdots w_k}\right]  f(\m{z})\, dy_1\cdots dy_{2k}.
\]
Now, we recall that $(\pi z)^{-1}$ is the fundamental solution of $\pa_{\ba{z}}$, see e.g. \cite{MR0203075}. Then we can substitute in the above formula $ \pa_{\ba{w}_j}\left[\frac{1}{w_j}\right] = \pi \, \del(w_j)=\pi \del(y_j) \del(y_{k+j})$. Finally, we obtain
\begin{align*}
\left\langle \frac{\del(\m{z})}{\det(A)}, f \right\rangle &=\frac{1}{\det(A)} \int_{D^k}  \del(\m{w})  f(\m{z})\, dy_1\cdots dy_{2k}\\
&= \ba{\det({A})} \int_{A^{-1} D^k}  \del(\m{w})  f(\m{z})\, dx_1\cdots dx_{2k}\\
&= \left\langle {\ba{\det({A})} \, \del(\m{w})}, f \right\rangle,
\end{align*}
which yields (\ref{delCompx}). In the second equality we have used the fact that $dy_1\cdots dy_{2k}= \det(\Psi(A))\, dx_1\cdots dx_{2k}= \det({A}) \ba{\det({A})} \,dx_1\cdots dx_{2k}$. 

\section*{Acknowledgements}
The authors want to thank Dr. Michael Wutzig for the careful reading of the manuscript and his valuable suggestions. Alí Guzmán Adán is supported by a BOF-post-doctoral grant from Ghent University.

\bibliographystyle{abbrv}

\end{document}